\documentclass[a4paper, 11pt]{article}
\usepackage{fullpage}
\usepackage{amsmath}
\usepackage{amssymb}
\usepackage{amsfonts}
\usepackage[french]{babel}
\usepackage{epsf}
\usepackage{epsfig}

\newtheorem{theo}{\indent Th\'eor\`eme\newline}[section]
\newtheorem{prop}[theo]{\indent Proposition\newline}
\newtheorem{lemme}[theo]{\indent Lemme\newline}
\newtheorem{cor}[theo]{\indent Corollaire\newline}

\def\N{{\Bbb{N}}}
\def\Z{{\Bbb{Z}}}
\def\R{{\Bbb{R}}}
\def\C{{\Bbb{C}}}
\def\cp{\C P^1}

\newcommand{\ind}{\mathop{\rm ind}\nolimits} 

\begin{document}

\title{Orientations complexes des $J$-courbes r\'eelles \\
Complex orientations of real $J$-curves}
\author{Jean-Yves Welschinger \\
Institut de Recherche Math\'ematique Avanc\'ee \\
Universit\'e Louis Pasteur \\
7, rue Ren\'e Descartes \\
67084 Strasbourg C\'edex, France \\
e-mail : welsch@math.u-strasbg.fr}
\maketitle

\makeatletter\renewcommand{\@makefnmark}{}\makeatother
\footnote{\hspace*{-20pt}Mots-cl\'es :  courbe alg\'ebrique r\'eelle, vari\'et\'e
presque-complexe.}  
\footnotetext{\hspace*{-20pt}Classification AMS : 14P25, 14Q05.}

{\bf R\'esum\'e :}

L'objet de cet article est d'obtenir des
restrictions sur les orientations complexes des $J$-courbes r\'eelles s\'eparantes
des vari\'et\'es presque-complexes. Ce probl\`eme provient de la
g\'eom\'etrie alg\'ebrique r\'eelle et le r\'esultat principal
est une congruence qui g\'en\'eralise des r\'esultats 
ant\'erieurs d'Arnol'd, Rokhlin, Mishachev, Zvonilov et Mikhalkin (v. 
\cite{Arn}, \cite{Rokh}, \cite{Mish}, \cite{Zvo}, \cite{Mikh}). Cette
congruence est obtenue \`a l'aide d'une application canonique d\'efinie sur les 
groupes d'homologies de dimension moiti\'e de ces vari\'et\'es \`a coefficients 
dans $\Z / l \Z$, \`a valeurs dans $\Z / 2l \Z$ (v. \S 1). \\

{\bf Abstract :}

The aim of this article is to obtain restrictions on complex orientations
of dividing real $J$-curves of almost-complex manifolds. This problem comes
from real algebraic geometry and the main result is a congruence generalising
Arnol'd, Rokhlin, Mishachev, Zvonilov and Mikhalkin 's previous results (see 
\cite{Arn}, \cite{Rokh}, \cite{Mish}, \cite{Zvo}, \cite{Mikh}). This
congruence is obtained with a canonical map defined on homology groups
of dimension half the dimension of the manifold with coefficients
in  $\Z / l \Z$, with values in $\Z / 2l \Z$ (see \S 1).

\section*{Introduction}

Soit $X$ une vari\'et\'e de classe $C^1$, de dimension $2n$, munie d'une structure 
presque-complexe $J$ de classe $C^1$. Une {\it structure r\'eelle} $c$ sur $X$ 
est une 
involution de classe $C^1$ dont la diff\'erentielle anti-commute avec
$J$. L'ensemble des points fixes de $c$ est alors une sous-vari\'et\'e de 
dimension $n$ de $X$ not\'ee $X_\R$, et appel\'ee
{\it partie r\'eelle} de $X$ (cette terminologie provient de la g\'eom\'etrie 
alg\'ebrique r\'eelle). Des exemples de telles vari\'et\'es sont donn\'es par $\C^k$, 
$\C P^k$ munis de la conjugaison complexe,
$\C P^k \times \C P^l$ muni du produit des conjugaisons complexes, ou encore
toutes leurs sous-vari\'et\'es complexes qui sont stables par la conjugaison
complexe.

Une {\it $J$-courbe r\'eelle} de $X$ est un plongement dans $X$ d'une surface de
Riemann $S$, o\`u la surface $S$ est munie d'une structure r\'eelle $c_S$, 
et le plongement commute d'une part avec
les involutions $c$ et $c_S$, et d'autre part sa diff\'erentielle commute avec 
les
structures presque-complexes de $S$ et $X$. Remarquons que la surface $S$
n'est pas forc\'ement connexe et que
ce plongement se restreint en un plongement de $S_\R$ dans $X_\R$. Une
premi\`ere question se pose alors : quels sont, \`a isotopie de $X_\R$ pr\`es, les
courbes de $X_\R$ obtenues comme parties r\'eelles de $J$-courbes r\'eelles
r\'ealisant une classe d'homologie de $H_2 (X ; \Z)$ donn\'ee ? 

Cette question a \'et\'e pos\'ee par D. Hilbert pour les courbes alg\'ebriques r\'eelles
non-singuli\`eres de $\C P^2$ dans son $16^{eme}$ probl\`eme (v. \cite{Hilb}), elle
est toujours ouverte \`a ce jour \`a partir du degr\'e $8$ (voir les surveys
\cite{Wils} et \cite{Viro}). Un premier r\'esultat est que le 
nombre de
composantes connexes de la partie r\'eelle $S_\R$ ne d\'epasse pas $g+\mu$ o\`u $g$
est le genre de la courbe $S$ (c'est-\`a-dire la somme des genres de ses
composantes connexes), et $\mu$ est le nombre de composantes connexes de $S$. 
Ce genre est bien d\'efini par la classe d'homologie
r\'ealis\'ee par la $J$-courbe, puisqu'une telle courbe v\'erifie la formule 
d'adjonction. Les
courbes poss\'edant ce nombre maximal de composantes sont appel\'ees {\it courbes
maximales}, ou $M$-courbes ; celles dont le nombre de composantes diff\`ere
du nombre maximal de $r$ sont not\'ees $(M-r)$-courbes.

Une surface de Riemann r\'eelle $S$ est dite {\it s\'eparante} s'il existe une
partie $S^+$ de $S$, bord\'ee par $S_\R$, telle que $S^+ \cup c_S (S^+) = S$
et $S^+ \cap c_S (S^+) = S_\R$. 
Une telle partie $S^+$ (qui est canoniquement 
orient\'ee, l'orientation provenant de la structure presque-complexe de $S$) 
induit une orientation sur $S_\R$ appel\'ee 
{\it orientation complexe} ; lorsque $S$ est connexe par exemple, sa partie 
r\'eelle $S_\R$ poss\`ede exactement deux orientations complexes, qui sont 
oppos\'ees.
Une $J$-courbe r\'eelle est dite {\it s\'eparante} si c'est un plongement d'une 
surface de Riemann s\'eparante, et une seconde question appara\^{\i}t : quels 
sont, \`a isotopie de $X_\R$
pr\`es, les plongements de courbes orient\'ees dans $X_\R$ obtenus comme parties 
r\'eelles de $J$-courbes r\'eelles s\'eparantes r\'ealisant une classe d'homologie 
de $H_2 (X ; \Z)$ donn\'ee ?

Les premi\`eres restrictions sur ces plongements ont \'et\'e donn\'ees par une 
congruence d'Arnol'd (\cite{Arn}), et une formule de Rokhlin (\cite{Rokh}).
Ces restrictions s'appliquent aux courbes
alg\'ebriques r\'eelles s\'eparantes non-singuli\`eres de degr\'es pairs de $\C P^2$ ;
toutefois, la formule de Rokhlin a \'et\'e \'etendue aux courbes planes 
projectives non-singuli\`eres de degr\'es impairs par Mishachev 
(\cite{Mish}). Plusieurs r\'esutats ont \'et\'e obtenus pour les courbes sur 
l'hyperbolo\"{\i}de (v. \cite{Gud}, \cite{Mat}, \cite{Mikh}), en particulier 
Mikhalkin a obtenu quelques
obstructions sous forme de congruences (v.\cite{Mikh} th\'eor\`emes\footnote{Il 
faut remplacer $dt+rs \equiv 0 \, \mod (4)$ par 
$dt+rs \equiv 2 \, \mod (4)$ dans le th\'eor\`eme 5.9.} 5.5, 5.7, 
5.8, 5.9 et 6.1 d). Le cas des surfaces complexes plus g\'en\'erales a \'et\'e
peu \'etudi\'e, Zvonilov (\cite{Zvo}) a g\'en\'eralis\'e une partie\footnote{Pour 
obtenir une formule d'orientations complexes g\'en\'eralisant celle de Rokhlin, 
il faut
calculer les classes d'homologies des cycles construits par Zvonilov, et leur 
indice d'intersection. Un tel calcul est effectu\'e dans le \S\ref{calcul}, mais
pour des surfaces $X$ plus particuli\`eres.} du r\'esultat de 
Rokhlin aux surfaces compactes complexes quelconques, mais pour les
courbes dont la partie r\'eelle r\'ealise la classe d'homologie nulle dans
$H_1 (X_\R ; \Z)$.

Dans cet article, tous ces r\'esultats sont unifi\'es et g\'en\'eralis\'es en une seule 
et m\^eme congruence (v. corollaire 2.6). Le \S$1$ contient quelques 
d\'efinitions, en particulier celle de l'application ${\cal P}_c$ qui est
fondamentale dans la suite. Soit $l$ un entier pair, la forme quadratique 
$x \mapsto x \circ c_* (x)$ de $H_n (X ; \Z)$ passe au quotient
en une application $q_c \, : \, H_n (X ; \Z) \otimes \Z /l\Z \to \Z /2l\Z$.
L'application ${\cal P}_c$ \'etend cette forme \`a $H_n (X ; \Z/l\Z)$, elle
est d\'efinie g\'eom\'etriquement. Les principaux r\'esultats de cet article, ainsi
que les liens avec les r\'esultats pr\'ec\'edents sont
\'enonc\'es et d\'emontr\'es dans le \S$2$. Le \S$3$ est consacr\'e au cas des surfaces 
fibr\'ees (sans fibres singuli\`eres) dont la base ou les fibres sont simplement 
connexes, c'est le cas par exemple des surfaces r\'egl\'ees. Enfin dans le \S$4$ 
sont donn\'es quelques exemples concrets de classes d'isotopies de courbes 
orient\'ees de $X_\R$ qui 
ne sont pas r\'ealis\'es comme parties r\'eelles de $J$-courbes r\'eelles s\'eparantes 
d'une classe donn\'ee, ce qui appara\^{\i}t comme une cons\'equence des 
r\'esultats pr\'ec\'edents.\\

{\bf Remerciement :}

Ce travail a \'et\'e effectu\'e sous la direction de V. Kharlamov. Je lui en
suis vivement reconnaissant.\\

{\bf Notations :}

$T_x M$ : espace tangent en $x$ de la vari\'et\'e $M$ de classe $C^1$.

$\chi (A)$ : caract\'eristique d'Euler de $A$.

$\partial A$ : bord de la cha\^{\i}ne $A$.

$[A]$ : classe d'homologie du cycle $A$.

$x \circ y$ : forme bilin\'eaire d'intersection de Poincar\'e.

$C_k (X ; R)$, $Z_k (X ; R)$, $B_k (X ; R)$, $H_k (X ; R)$ :
espace des cha\^{\i}nes (resp. cycles, bords, classes d'homologies) singuli\`eres
de dimension $k$ de $X$ \`a coefficients dans $R$.

{\bf Toutes les cha\^{\i}nes seront enti\`eres, les cycles quant \`a eux seront 
entiers ou \`a coefficients dans $\Z / l\Z$.}

\section{Application ${\cal P}_c$}

Soit $X$ une vari\'et\'e de classe $C^1$ compacte orient\'ee de dimension $2n$ 
($n \in \N$),
 munie d'une involution $c$ de classe
$C^1$ qui pr\'eserve l'orientation et telle que chaque composante 
connexe de la sous-vari\'et\'e $X^c$ de ses points fixes 
est de dimension au plus $n$. 

\subsection{Cha\^{\i}nes simpliciales}

Une $k$-cha\^{\i}ne singuli\`ere $A$ de $X$ ($k \in \{ 0 , \dots , 2n \}$) est
dite {\it simpliciale} s'il existe une triangulation de classe $C^1$ de $X$
pour laquelle $A$ soit une cha\^{\i}ne simpliciale. Dans ce cas, le 
{\it squelette de dimension q} de $A$ associ\'e \`a cette triangulation 
est la r\'eunion des faces de $A$ 
de dimensions inf\'erieures ou \'egales \`a $q$ ($q \in \{ 0 , \dots , k \}$).

Deux cha\^{\i}nes simpliciales $A$ et $B$ sont dites {\it transverses} si en 
tout
point d'intersection $x$ de $A$ et $B$, les espaces tangents $T_A$ et $T_B$ 
des faces de $A$  et $B$ (associ\'ees \`a des triangulations de $X$)
 contenant $x$ sont en position g\'en\'erale dans $T_x X$. Consid\'erons 
deux suppl\'ementaires orient\'es $\omega_A$ et $\omega_B$ de $T_A$ et $T_B$ : 
c'est-\`a-dire deux sous-espaces vectoriels orient\'es de $T_x X$ tels que
$T_A \oplus \omega_A = T_B \oplus \omega_B = T_x X$ en tant qu'espaces 
orient\'es. L'intersection $T_A \cap T_B$ est orient\'ee de sorte
que $(T_A \cap T_B) \oplus \omega_A \oplus \omega_B = T_x X$ en tant 
qu'espaces orient\'es. Si $A$ et 
$B$ sont deux cha\^{\i}nes simpliciales transverses, il existe une 
triangulation de classe $C^1$ de $X$ pour laquelle l'intersection 
alg\'ebrique $A \circ B$ est \'egalement une cha\^{\i}ne simpliciale dont les 
faces sont les intersections orient\'ees des faces de $A$ et $B$,
et les multiplicit\'es de ces faces sont les produits des multiplicit\'es des 
faces de $A$ et $B$.

Une cha\^{\i}ne simpliciale $A$ de dimension $k$ d'une vari\'et\'e \`a bord $M$ est 
dite
{\it transverse au bord} de $M$ si dans un voisinage du bord de $M$, 
identifi\'e \`a un produit de $\partial M$ par un intervalle, la cha\^{\i}ne
$A$ s'identifie elle aussi \`a un produit de $A \cap \partial M$ par un 
intervalle.

Les deux lemmes suivants seront utiles dans la suite :

\begin{lemme}
\label{lemmebbord}
Soient $M$ une vari\'et\'e orient\'ee \`a bord de dimension $2n+1$ et $A$, $B$ deux 
cha\^{\i}nes simpliciales de dimension $n+1$ de $C_k (M ; R)$.
Supposons que $A$ et $B$ soient transverses au bord et que
$\partial A$ et $\partial B$ soient deux cha\^{\i}nes simpliciales
 transverses de $\partial M$. Alors en tout point 
$x$ de $\partial A \cap \partial B$ :
$$\partial (A \circ_M B)(x) = \partial A \circ_{\partial M} \partial B (x). 
\; \square$$
\end{lemme}

\begin{lemme}
\label{lemmebord}
Si $A$ et $B$ sont deux cha\^{\i}nes simpliciales transverses de dimension 
$n+1$ 
d'une vari\'et\'e $M$ orient\'ee de dimension $2n+1$, alors en tout point $x$ de 
$A \cap \partial B$,
$$\partial (A \circ B) (x) = A \circ \partial B (x). \; \square$$
\end{lemme}
(Dans le lemme $1.2$, la vari\'et\'e $M$ peut \'eventuelement \^etre \`a bord. Dans ce 
cas, la condition de transversalit\'e de $A$ et $B$ implique que le point $x$ 
est dans l'int\'erieur de $M$.) 

\subsection{Application ${\cal P}_c$}

 Soient $l$ un entier pair et $A$ une cha\^{\i}ne simpliciale de 
dimension $n$ de $X$, dont le bord est nul modulo $l$. Si la cha\^{\i}ne $A$
est telle que $A$ et $c(A)$ sont transverses, elle est dite 
{\it g\'en\'erique}. Pour de telles cha\^{\i}nes, on pose  
${\cal P}_c (A) = A \circ c(A)$ r\'eduit modulo $2l$ (de sorte que ${\cal P}_c$
est \`a valeur dans $\Z /2l\Z$).

\begin{prop}
\label{P_c}
Il existe des cycles g\'en\'eriques dans chaque classe d'homologie de 
$H_n (X ; \Z / l\Z )$, et ${\cal P}_c$ passe au quotient en une application 
$H_n (X ; \Z / l\Z ) \to \Z / 2l\Z $ (toujours not\'ee ${\cal P}_c$) telle
que le diagramme
$$\vcenter{\hbox{\begin{picture}(0,0)%
\epsfig{file=reg2.pstex}%
\end{picture}%
\setlength{\unitlength}{0.00083300in}%
\begingroup\makeatletter\ifx\SetFigFont\undefined%
\gdef\SetFigFont#1#2#3#4#5{%
  \reset@font\fontsize{#1}{#2pt}%
  \fontfamily{#3}\fontseries{#4}\fontshape{#5}%
  \selectfont}%
\fi\endgroup%
\begin{picture}(2337,1095)(3001,-4600)
\put(3001,-4561){\makebox(0,0)[lb]{\smash{\SetFigFont{12}{14.4}{\rmdefault}{\mddefault}{\updefault}$H_n (X ; \Z / l\Z)$}}}
\put(5101,-4561){\makebox(0,0)[lb]{\smash{\SetFigFont{12}{14.4}{\rmdefault}{\mddefault}{\updefault}$\Z / 2l\Z$}}}
\put(3526,-4111){\makebox(0,0)[lb]{\smash{\SetFigFont{12}{14.4}{\rmdefault}{\mddefault}{\updefault}$red_l$}}}
\put(4501,-4411){\makebox(0,0)[lb]{\smash{\SetFigFont{12}{14.4}{\rmdefault}{\mddefault}{\updefault}${\cal P}_c$}}}
\put(3001,-3661){\makebox(0,0)[lb]{\smash{\SetFigFont{12}{14.4}{\rmdefault}{\mddefault}{\updefault}$H_n (X ; \Z) \otimes \Z/l\Z$}}}
\put(4801,-3886){\makebox(0,0)[lb]{\smash{\SetFigFont{12}{14.4}{\rmdefault}{\mddefault}{\updefault}$q_c$}}}
\end{picture}
}}$$
soit commutatif.
\end{prop}
(Rappelons que la forme quadratique $(x,y) \mapsto x \circ c_* (y)$ r\'eduite
 modulo $2l$ et d\'efinie sur $H_n (X ; \Z)$ passe au quotient en une application
de $H_n (X ; \Z) \otimes \Z /l\Z$ ; cette application est not\'ee $q_c$.
Par ailleurs, convenons que $\Z /0\Z = \Z$).\\

{\bf D\'emonstration :} 

Dans chaque classe d'homologie de $H_n (X ; \Z / l\Z )$, un cycle simplicial
 qui est en position
g\'en\'erale dans $X$ est g\'en\'erique puisque par hypoth\`ese les
composantes connexes de $X^c$ sont de dimension au plus $n$.

Soient $A_0$, $A_1$ deux cha\^{\i}nes g\'en\'eriques de dimension $n$,
de bord nul modulo $l$ et qui sont homologues dans $C_n (X ; \Z / l\Z )$.
Il s'agit de montrer que ${\cal P}_c (A_0)$ et ${\cal P}_c (A_1)$ sont \'egaux.
Notons $\widetilde{X} = X \times [0,1]$,
$\tilde{c} = c\times id$, $\widetilde{X}^c = X^c \times [0,1]$, et plongeons
$A_i$ dans $X \times \{ i \}$ ($i \in \{0,1\}$).

Il existe une cha\^{\i}ne simpliciale $H \in C_{n+1} (\widetilde{X} ; \Z)$
telle que $\partial H = A_1 - A_0 - l \theta$ o\`u 
$\theta \in C_n (\widetilde{X} ; \Z)$ (et $\theta$ est \'egalement simplicial). 
La cha\^{\i}ne 
$H$ peut \^etre choisie en position g\'en\'erale dans $\widetilde{X}$ de sorte 
que :

\begin{itemize}
\item le squelette de dimension $n-1$ de $H$ n'intersecte pas $c(H)$.
\item $H$ est transverse \`a $\widetilde{X}^c$ et au bord de $\widetilde{X}$.
\item si $x$ est un point d'intersection de $H$ et $c(H)$ n'appartenant pas \`a
$\widetilde{X}^c$, alors $H$ et $c(H)$ sont transverses en $x$. 
\item si $x$ est un point d'intersection d'un simplexe de dimension $n+1$ de 
$H$ et de $\widetilde{X}^c$, alors $H$ et $c(H)$ sont transverses en $x$ ; et 
si $x$ est un point d'intersection d'un simplexe de dimension $n$ de 
$H$ et de $\widetilde{X}^c$,
l'espace tangent en $x$ de ce simplexe 
et de son conjugu\'e engendrent un sous-espace $E_x$ de dimension $2n$ de 
$T_x \widetilde{X}$ qui ne contient aucun espace tangent en $x$ de simplexes
 de dimension $n+1$ de $H$ contenant $x$.
\end{itemize}

L'intersection alg\'ebrique $H \circ c(H)$ est une $1$-cha\^{\i}ne \`a 
coefficients entiers,
de sorte que $\partial (H \circ c(H))$ est une $0$-cha\^{\i}ne dont la somme
des coefficients est nulle. Par
ailleurs, si $x$ est un point du bord de cette $1$-cha\^{\i}ne, soit
$x$ appartient \`a l'intersection de $A_1$ ou $A_0$ avec son conjugu\'e, soit $x$
est un point d'intersection de $l\theta$ avec un
simplexe de dimension $n+1$ de 
$c(H)$ ($x$ n'appartient alors pas \`a $\widetilde{X}^c$), soit $x$ est un point
d'intersection de $l\theta$ et $c(l\theta )$ ($x$ appartient 
alors \`a $\widetilde{X}^c$).
D'apr\`es le lemme \ref{lemmebbord}, la contribution au bord de $H \circ c(H)$ 
des points de la premi\`ere cat\'egorie est ${\cal P}_c (A_1) - {\cal P}_c (A_0)$
(l'orientation de $X \times \{ 0 \}$ est oppos\'ee \`a celle de $X$). D'apr\`es
le lemme \ref{lemmebord}, la contribution modulo $2l$ des points de la 
seconde cat\'egorie
est nulle (ces points sont coupl\'es par deux par l'involution $c$). Quant aux 
points de la troisi\`eme cat\'egorie, leur contribution
est \'egalement nulle modulo $2l$, en effet, soit $x$ un tel 
point d'intersection
et $E_x$ le sous-espace de dimension $2n$ de $T_x \widetilde{X}$ engendr\'e 
par les espaces tangents $T_1$ et $T_2$ en $x$ des simplexes de dimension 
$n$ de $H$ et $c(H)$.
Ce sous-espace $E_x$ s\'epare $T_x \widetilde{X}$
en deux demi-espaces $\widetilde{X}_1$ et $\widetilde{X}_2$ et
on note $H_1 = T_x H \cap \widetilde{X}_1$ et $H_2 = T_x H \cap 
\widetilde{X}_2$. Il existe deux entiers $t_1 , t_2$ tels que
 $\partial H_1 = t_1 T_1$, $\partial H_2 = t_2 T_1$ et $t_1 +t_2 =
0 \quad \mod (l)$. En appliquant le lemme \ref{lemmebbord}, on s'aper\c{c}oit que 
la contribution
au bord de $H \circ c(H)$ de $x$ est $(t_1^2 - t_2^2) T_1
 \circ_{\partial \widetilde{X}_1} T_2$, cette contribution est donc nulle 
modulo $2l$. $\square$\\

Les propri\'et\'es suivantes de l'application ${\cal P}_c$ se v\'erifient 
facilement :
\begin{prop}
  \begin{enumerate}
  \item ${\cal P}_c (a+b) = 
{\cal P}_c (a) + {\cal P}_c (b) + i(a \circ c_* (b))$ o\`u $i$ est l'injection
$\Z /l\Z \to \Z /2l\Z$ et $a,b \in H_n (X ; \Z /l\Z)$.
  \item ${\cal P}_c (\rho (a)) = a \circ c_* (a)$ o\`u $a \in H_n 
(X ; \Z /2l\Z)$ et $\rho$ : $H_n (X ; \Z /2l\Z) \to
H_n (X ; \Z /l\Z)$ est la r\'eduction modulo $l$.
\item Si $a \in H_n (X ; \Z /l\Z)$, alors ${\cal P}_c (a) = 
{\cal P}_c (c_*(a))$ et $a \circ c_* (a) = {\cal P}_c (a) \; \mod (l)$. 
$\square$
  \end{enumerate}
\end{prop}

\section{Orientations complexes des $J$-courbes r\'eelles}

\subsection{R\'esultat principal : une congruence}

Soit $X$ une vari\'et\'e presque-complexe compacte de dimension $4$ (r\'eelle), 
munie d'une structure r\'eelle $c$ et de
l'orientation induite par la structure presque-complexe.

Soient $\Lambda$ une $J$-courbe r\'eelle s\'eparante de $X$,
et $\Lambda^+$ une partie de $\Lambda \setminus \Lambda_\R$ bord\'ee par 
$\Lambda_\R$, qui induit une
orientation complexe sur $\Lambda_\R$ que nous fixons. Choisissons $j$ 
courbes lisses
orient\'ees $\lambda_1$, $\ldots$, 
$\lambda_j$ de $X_\R$, chacune soit disjointe
de $\Lambda_\R$, soit confondue avec une composante de $\Lambda_\R$, ainsi que
 $j+1$
entiers $l, s_1, \dots, s_j$ tels que $[\Lambda_\R] = l\sum_{i=1}^j s_i 
[\lambda_i] \in H_1 (X_\R ; \Z)$. Les courbes $\lambda_i$ peuvent \^etre choisies
de sorte que lorsqu'une telle courbe est confondue avec une composante de 
$\Lambda_\R$, les orientations de la courbe et de la composante 
co\"{\i}ncident.

Soit $\{ B_k \, , \, k \in K\}$ l'ensemble des composantes connexes du 
compl\'ementaire de
$\Lambda_\R \cup \cup_{i=1}^j \lambda_i $ dans $X_\R$. Lorsqu'une composante
$B_k$ est orientable, elle est munie d'une orientation arbitraire. De
la suite exacte longue associ\'ee \`a la paire $(X_\R , \Lambda_\R \cup 
\cup_{i=1}^k \lambda_i )$ d\'ecoule le lemme suivant :

\begin{lemme}
\label{tk}
Il existe une famille d'entiers $(t_k)_{k \in K}$ telle que $t_k = 0$ 
si $B_k$ n'est pas orientable, et 
$\partial (\sum_{k \in K} t_k B_k) = l \sum_{i=1}^j s_i \lambda_i - 
\Lambda_\R$. $\square$

\end{lemme}

Soient $( t_k)_{k \in K}$ une famille d'entiers donn\'ee par le lemme \ref{tk},
et $\ind$ la fonction
de $X_\R$ d\'efinie par 
$\ind(x) = t_k $ si $x \in B_k$ et $\ind(x) = 
0$ si $ x \in \Lambda_\R \cup \cup_{i=1}^j \lambda_i$. La somme
$\sum_{k \in K} t_k^2 \chi(B_k)$ est not\'ee 
$\int_{X_\R} \ind^2 d\chi$ (notation introduite par Viro dans \cite{not}).

\begin{theo}
\label{cong}
Soit $\Lambda$ une $J$-courbe r\'eelle s\'eparante de $X$,
telle que $[\Lambda_\R] = l\sum_{i=1}^j s_i [\lambda_i] \in H_1 (X_\R ; \Z)$
o\`u $l, s_1, \dots, s_j$ sont $j+1$ entiers et $\lambda_1, \ldots, \lambda_j$
sont $j$ courbes lisses orient\'ees de $X_\R$, chacune soit disjointe
de $\Lambda_\R$, soit incluse dans $\Lambda_\R$ auquel cas l'orientation de 
la courbe co\"{\i}ncide avec l'orientation de $\Lambda_\R$. Soient 
$( t_k)_{k \in K}$
une famille d'entiers donn\'ee par le lemme 2.1, $\ind$ la fonction de $X_\R$ 
associ\'ee et ${\cal A}_l$ le $2$-cycle $\Lambda^+ + \sum_{k \in K} t_k B_k
\in Z_2 (X ; \Z / l\Z)$. Si $l$ est pair, alors :
$${\cal P}_c({\cal A}_l) = - \int_{X_\R} \ind^2 d\chi \quad \mod (2l).$$
Si $l$ est impair, 
$[{\cal A}_l] \circ c_* [{\cal A}_l ] = - \int_{X_\R} \ind^2 d\chi 
\quad \mod (l).$
\end{theo}

{\bf D\'emonstration :}

Choisissons un champ de vecteurs $\xi$ de classe $C^1$, tangent \`a 
$\Lambda_\R \cup \cup_{i=1}^j \lambda_i$, compatible avec l'orientation de 
cette courbe. On \'etend ce champ \`a 
${\cal A}_l$ tout entier de sorte qu'en d\'ecalant ce cycle \`a l'aide d'un 
flot 
associ\'e au champ de vecteurs $J(\xi)$, on obtient un cycle simplicial
$\widetilde{\cal A}_l$ g\'en\'erique. Da la m\^eme fa\c{c}on que dans \cite{Rokh}, on 
s'aper\c{c}oit que $\widetilde{\cal A}_l \circ c ( \widetilde{\cal A}_l )
= - \int_{X_\R} \ind^2 d\chi$. $\square$

\begin{theo}
Soit $\Lambda$ une $J$-courbe r\'eelle s\'eparante  de $X$,
telle que $[\Lambda_\R] = 0 \in H_1 (X_\R ; \Z /2\Z)$. Soit
$B$ une partie 
$X_\R \setminus \Lambda_\R$ telle que $\partial B = \Lambda_\R \quad \mod (2)$
et ${\cal A}_l$ le cycle $\Lambda^+ + B \in Z_2 (X ;\Z /2\Z)$. Alors :
$$ {\cal P}_c({\cal A}_l) = -\chi(B) \quad \mod (4).$$
\end{theo}

{\bf D\'emonstration :}

Si la partie $B$ est orientable, c'est un corollaire du th\'eor\`eme \ref{cong}
 ; sinon,
ce th\'eor\`eme se d\'emontre comme le pr\'ec\'edent. $\square$\\

{\bf Remarque :} Si $Tor_1 (H_1 (X ; \Z) , \Z / l\Z) = 0$, le cycle
${\cal A}_l$ se rel\`eve en un cycle entier et d\'efinit un \'el\'ement $a$ de
$H_2 (X ; \Z) \otimes \Z / l\Z$. Dans ce cas la proposition \ref{P_c} implique
que ${\cal P}_c({\cal A}_l) = q_c (a)$.

\subsection{Calcul de la classe $[{\cal A}_l]$}
\label{calcul}
L'objet de ce paragraphe est de donner une classe de vari\'et\'es
presque-complexes pour lesquelles il est possible de calculer $a$ (et 
$q_c (a)$) en fonction de $[\Lambda] \in H_2 (X ; \Z)$, des orientations 
complexes de $\Lambda_\R$ et du choix des entiers $(t_k)_{k \in K}$.

\subsubsection{L'indice $i (\lambda ,{\cal A}_l )$}

Soient $\lambda$ une courbe lisse de $X_\R$, transverse \`a $\Lambda_\R$ 
et ${\cal A}_l$ le cycle d\'efini dans le th\'eor\`eme 2.2. Soit $\xi$ une section
de $(T X_\R)\left| \right. _{\lambda}$, de classe $C^1$, ne s'annulant pas, 
et transverse \`a 
$T \lambda \subset (T X_\R)\left| \right. _{\lambda}$. Supposons
de plus que les points en lesquels le champ $\xi$ est tangent \`a la courbe 
$\lambda$ n'appartiennent pas \`a 
$\Lambda_\R \cup \cup_{i=1}^j \lambda_i$, et notons $\Lambda_\R '$ un champ de 
vecteurs tangents \`a $\Lambda_\R$ ne s'annulant pas et compatible avec 
l'orientation de cette courbe.

Dans les sections 2.1 et 2.2, la courbe $\Lambda^+$ intervient dans le cycle
${\cal A}_l$ avec une multiplicit\'e $1$, tandis
que dans la section 2.3, cette courbe intervient dans ${\cal A}_l$ avec une 
multiplicit\'e sup\'erieure, not\'ee $q \in \N^*$. Soit $x$ un 
point d'intersection de $\lambda$ et $\Lambda_\R$, posons :
$$i_\xi (\lambda ,{\cal A}_l ) (x) = 
\left\{
  \begin{array}{l}
+q \quad \mbox{si $\xi (x)$ et $\Lambda_\R ' (x)$ appartiennent \`a la m\^eme
composante de $T_x X_\R \setminus T_x \lambda$}, \\
-q \quad \mbox{sinon.}
  \end{array}
\right.$$

Soit \`a pr\'esent $x$ un point en lequel $\xi (x)$ est tangent \`a $\lambda$, et 
$S$ une 
sph\`ere de $X_\R$ centr\'ee en $x$, choisie suffisament petite de sorte 
d'une part qu'elle 
soit incluse dans une partie $B_k$ et qu'elle ne contienne aucun 
autre point en lequel $\xi$ est tangent \`a
$\lambda$, et d'autre part que $\lambda$ ne l'intersecte qu'en deux points 
$x_1$ et $x_2$. Choisissons \'egalement $S$ de sorte que $\xi (x_1)$ et 
$\xi (x_2)$ soient tangents \`a $S$. Si $B_k$
n'est pas orientable, on pose $i_\xi (\lambda ,{\cal A}_l ) (x) = 0$. Sinon,
la sph\`ere $S$ est le bord d'un disque orient\'e de $B_k$ 
centr\'e en $x$ et h\'erite d'une orientation provenant de ce disque. Posons :
\begin{itemize}
\item $i_\xi (\lambda ,{\cal A}_l ) (x) = 0$ si $\xi (x_1)$ et $\xi (x_2)$ 
ne se prolongent pas en un champ de
vecteurs tangents de $S$.
\item $i_\xi (\lambda ,{\cal A}_l ) (x) = 2t_k$ si $\xi (x_1)$ et $\xi (x_2)$ 
se prolongent en un champ de
vecteurs tangents de $S$ compatible avec son orientation, et 
$i_\xi (\lambda ,{\cal A}_l ) (x) = -2t_k$ dans le cas contraire.
\end{itemize}

Finalement, posons $i_\xi (\lambda ,{\cal A}_l ) = \sum
i_\xi (\lambda ,{\cal A}_l ) (x)$ (la somme \'etant prise sur les points $x$
d'intersection de $\lambda$ et $\Lambda_\R$ et les points en lesquels $\xi$
est tangent \`a $\lambda$). Cette valeur ne change pas si on
modifie $\xi$ en dehors de $\Lambda_\R \cup \cup_{i=1}^j \lambda_i$.

\begin{lemme}
\label{ixi}
La valeur modulo $2l$ de $i_\xi (\lambda , {\cal A}_l )$ est ind\'ependante du 
champ $\xi$. Par ailleurs, si le champ $\xi$ est fix\'e sur 
$(\cup_{i=1}^j \lambda_i) \cap \lambda$, 
alors $i_\xi (\lambda ,{\cal A}_l )$ est un entier qui ne d\'epend pas de la 
fa\c{c}on de 
prolonger ce champ \`a $\lambda$.
\end{lemme}

On note $i (\lambda ,{\cal A}_l ) \in \Z/2l\Z $ l'\'el\'ement donn\'e par le lemme
\ref{ixi}.

{\bf D\'emonstration :}

Choisissons une famille de champs de vecteurs de r\'ef\'erence de classe $C^1$, 
d\'efinis sur $\lambda$, tous \'egaux en dehors d'un voisinage d'un point $x_0$ 
fix\'e (n'appartenant pas \`a $\Lambda_\R \cup \cup_{i=1}^j \lambda_i$), ne 
s'annulant pas, et qui ne sont tangents \`a $\lambda$ qu'au voisinage de ce 
point $x_0$.

Pour v\'erifier la premi\`ere partie du lemme, on s'aper\c{c}oit que pour tout champ de
vecteurs $\xi$ de classe $C^1$ ne s'annulant pas et transverse \`a $T \lambda
\subset (T X_\R)\left| \right. _{\lambda}$, il existe un chemin continu 
$\xi (t)$ de
champs de vecteurs de classe $C^1$ ne s'annulant pas et transverse \`a $T \lambda
\subset (T X_\R)\left| \right. _{\lambda}$ reliant $\xi$ \`a l'un des champs 
de vecteurs de r\'ef\'erence, et que
la valeur modulo $2l$ de $i_{\xi (t)} (\lambda ,{\cal A}_l )$ ne d\'epend pas 
de $t$ (lorsqu'elle est bien d\'efini). Or la valeur de $i_\xi 
(\lambda ,{\cal A}_l )$ pour un
champ qui est d\'efini en dehors d'un voisinage de $x_0$ ne d\'epend pas du 
prolongement de ce champ \`a $\lambda$ tout entier. La seconde partie du
lemme \ref{ixi} se v\'erifie de fa\c{c}on analogue. $\square$

\subsubsection{Calcul de la classe $[{\cal A}_l]$}

Soit $X$ une vari\'et\'e presque-complexe compacte de dimension $4$ 
munie d'une structure r\'eelle $c$ et de l'orientation induite par la structure 
presque-complexe. Choisissons des $J$-courbes r\'eelles $v^1, \dots , v^n$ dans
 $X$.
Supposons que la $J$-courbe $\Lambda$
donn\'ee par le th\'eor\`eme \ref{cong} soit transverse \`a $v^1, \dots , v^n$, et 
notons 
$i ({\cal A}_l )$ l'\'el\'ement de $(\Z/l\Z)^n$ dont les coordonn\'ees sont 
$\frac{1}{2} (i (v^m_\R , {\cal A}_l )+ v^m \circ \Lambda)$ ($m \in \{ 1 , 
\dots ,n \}$). Notons \'egalement $Q$ la matrice de la forme d'intersection 
de Poincar\'e de $X$ dans la famille $([v^1], \dots , [v^n])$.

\begin{theo}
\label{homol}
Supposons que
$Tor_1 (H_1 (X ; \Z) ; \Z / l\Z) = 0$ et que $([v^1], \dots , [v^n])$ forme une
base de $H_2 (X ; \Z)/Tors H_2 (X ; \Z)$. Soit $\Lambda$ une $J$-courbe 
r\'eelle de $X$ satisfaisant aux hypoth\`eses du th\'eor\`eme \ref{cong} et transverse
aux courbes $v^1, \dots , v^n$.
En utilisant les notations du th\'eor\`eme \ref{cong}, on a :
$$[{\cal A}_l] =  [Q^{-1} 
(i ({\cal A}_l ))]  \in H_2 (X ; \Z) \otimes \Z / l\Z $$
\end{theo}
(Ici, $[Q^{-1} (i ({\cal A}_l ))]$  d\'esigne la classe de $H_2 (X ; \Z) 
\otimes \Z / l\Z$ dont les coordonn\'ees dans la base $([v^1], \dots , [v^n])$
sont $Q^{-1} (i ({\cal A}_l ))$.)

Les th\'eor\`emes \ref{cong}, 2.3, \ref{homol} et la proposition \ref{P_c} 
entra\^{\i}nent les corollaires suivants :

\begin{cor}

Sous les hypoth\`eses et notations du th\'eor\`eme 2.5,
$$\int_{X_\R} \ind^2 d\chi = -q_c \big( [Q^{-1} 
(i ({\cal A}_l ))]  \big) \quad \left[ 
\begin{array}{l} 
\mod (2l) \mbox{ si $l$ est pair,} \\
\mod (l) \mbox{ sinon. }\square
\end{array} \right.$$
\end{cor}

{\bf Remarque :} Si $l$ est impair, pour d\'eterminer l'int\'egrale
$\int_{X_\R} \ind^2 d\chi$ modulo $2l$ il suffit de la d\'eterminer modulo $l$,
puisque sa valeur modulo $2$ est connue (c'est la caract\'eristique d'Euler de 
la r\'eunion des parties $B_k$ pour lesquelles $t_k$ est impair). 

\begin{cor}
Supposons que $X$ satisfasse aux hypoth\`eses du th\'eor\`eme 2.5 et soit $\Lambda$ 
une $J$-courbe r\'eelle s\'eparante 
de $X$ telle que $[\Lambda_\R] = 0 \in H_1 (X_\R ; \Z /2\Z)$. 
Soit $B$ une partie de 
$X_\R \setminus \Lambda_\R$ telle que $\partial B = \Lambda_\R \quad 
\mod (2)$, et ${\cal A}_l$ le cycle $\Lambda^+ + B \in Z_2 (X ;\Z /2\Z)$. 
Alors :
$$ \chi(B) = -q_c \big([Q^{-1} (i ({\cal A}_l ))] \big) \quad \mod (4). 
\; \square$$
\end{cor}

{\bf D\'emonstration du th\'eor\`eme 2.5 :}

Rappellons que $\partial {\cal A}_l = l\sum_{i=1}^j s_i \lambda_i$. Soit $\xi$
un champ de vecteurs de classe $C^1$ de $X_\R$ ne s'annulant pas, d\'efini sur 
$\cup_{i=1}^j \lambda_i \cup \cup_{m=1}^n v^m_\R$, et qui n'est tangent aux 
courbes $v^m$ qu'en des points n'appartenant pas \`a 
$\Lambda_\R \cup \cup_{i=1}^j \lambda_i$.
Les courbes $\lambda_i$ peuvent \^etre d\'ecal\'ees \`a l'aide d'un flot associ\'e au 
champ de vecteurs $- J( \xi)$ de fa\c{c}on \`a
obtenir des courbes $\tilde{\lambda}_i$ disjointes de $X_\R$.

La classe d'homologie $\sum_{i=1}^j s_i [\tilde{\lambda}_i] \in H_1 (X ; \Z)$ 
est d'ordre $l$, et $Tor_1 (H_1 (X ; \Z) ; \Z / l\Z) =0$ par hypoth\`ese.
Il existe donc un $2$-cycle $\widetilde{S}$ de $X$ tel que 
$\partial \widetilde{S} =
\sum_{i=1}^j s_i \tilde{\lambda}_i$, et on prolonge ce cycle \`a l'aide du flot
pr\'ec\'edent en un cycle $S$ tel que
${\cal A} =  {\cal A}_l - lS$ soit un $2$-cycle entier de $X$ qui rel\`eve
${\cal A}_l$.

Consid\'erons \`a pr\'esent une extension du champ $\xi$ (toujours not\'ee $\xi$)
aux courbes $v^m$ ($m \in \{ 1 , 
\dots ,n \}$), et d\'ecalons ces courbes \`a l'aide d'un flot associ\'e au champ 
de vecteurs $J(\xi) $. On obtient ainsi $n$ courbes $\tilde{v}^m$ 
($m \in \{ 1 , \dots ,n \}$), 
homologues aux courbes $v^m$ et dont les intersections avec ${\cal A}$ sont
transverses. Soit $\tilde{x}$ un point d'intersection de $\tilde{v}^m$ et
${\cal A}$, ce point
correspond soit \`a un point d'intersection $x$ de $v^m$ avec $\Lambda$, soit
\`a un point $x$ en lequel $\xi $ est tangent \`a $v^m_\R$, soit \`a un point
d'intersection $x$ de $v^m$ avec $S$, et :
\begin{itemize}
\item si $\tilde{x}$ appartient \`a $ \Lambda$ et n'est pas voisin de $X_\R$, 
alors
$\tilde{v}^m \circ {\cal A} (\tilde{x}) = v^m \circ \Lambda (x) =
\frac{1}{2} (v^m \circ \Lambda (x) + v^m \circ \Lambda (c(x)))$.
\item si $\tilde{x}$ est voisin de $X_\R$ et appartient \`a $ \Lambda$, 
alors $\tilde{v}^m \circ {\cal A} (\tilde{x}) = v^m \circ \Lambda (x) \, 
\mod (l) =
\frac{1}{2} (v^m \circ \Lambda (x) + i_\xi (v^m_\R , {\cal A}_l) (x))
\, \mod (l)$ (ce sont des \'egalit\'es entre entiers si le point $x$ 
n'appartient pas en plus \`a une courbe $\lambda_i$).
\item si $\tilde{x}$ appartient \`a $S$ et correspond \`a un point $x$ 
n'appartenant pas \`a $\Lambda$, alors $\tilde{v}^m \circ {\cal A} 
(\tilde{x}) = 0 \quad \mod (l)$.
\item si $\tilde{x}$ n'appartient ni \`a $ \Lambda$, ni \`a $S$, alors
$\tilde{v}^m \circ {\cal A} (\tilde{x}) =\frac{1}{2} i_\xi (v^m_\R , 
{\cal A}_l) (x)$.
\end{itemize}

Il en d\'ecoule que $\tilde{v}^m \circ {\cal A} = \frac{1}{2} 
(v^m \circ \Lambda  + i_\xi (v^m_\R , {\cal A}_l) ) \,
\mod (l)$ et donc que
$[{\cal A}_l] = [Q^{-1} 
(i ({\cal A}_l ))] \in H_2 (X ; \Z) \otimes \Z / l\Z$,
puisque $([v^1], \dots , [v^n])$ est une base de $H_2 (X ; \Z) \otimes \Z / 
l\Z$. $\square$

\subsection{Compl\'ements.}

\subsubsection{Compl\'ement au corollaire $2.6$.}

Soit $X$ une vari\'et\'e presque-complexe compacte de dimension $4$ (r\'eelle) 
satisfaisant 
aux hypoth\`eses du th\'eor\`eme 2.5, et $\Lambda$ une $J$-courbe r\'eelle s\'eparante 
de $X$. La courbe $\Lambda_\R$ est munie d'une orientation complexe
provenant d'une partie $\Lambda^+$ de $\Lambda \setminus \Lambda_\R$ bord\'ee
par $\Lambda_\R$.

Supposons qu'il existe $j+2$ entiers $q$, $l, s_1, \dots, s_j$, et
$j$ courbes lisses orient\'ees $\lambda_1$, $\ldots$, 
$\lambda_j$ de $X_\R$, chacune soit disjointe
de $\Lambda_\R$, soit confondue avec une composante de $\Lambda_\R$, tels que
$q [\Lambda_\R] = q l \sum_{i=1}^j s_i [\lambda_i] \in 
H_1 (X_\R ; \Z)$. Les courbes $\lambda_i$ peuvent \^etre choisies
de sorte que lorsqu'une telle courbe est confondue avec une composante de 
$\Lambda_\R$, les orientations de la courbe et de la composante 
co\"{\i}ncident.

Notons $\{ B_k \, , \, k \in K\}$ l'ensemble des composantes connexes du 
compl\'ementaire de $\Lambda_\R \cup \cup_{i=1}^j \lambda_i $ dans $X_\R$ ; 
il existe une famille $(t_k)_{k \in K}$ d'entiers
tels que $\partial (\sum_{k \in K} t_k B_k) = q l \sum_{i=1}^j s_i 
\lambda_i - q \Lambda_\R$ (de la m\^eme fa\c{c}on que dans le lemme $2.1$). Soient 
$\ind$ la fonction de $X_\R$ associ\'ee \`a 
cette famille (v. \S2.1 pour une d\'efinition), et ${\cal A}_{q , l}$ le cycle 
$q \Lambda^+ + \sum_{k \in K} t_k B_k \in Z_2 (X ; \Z / q l \Z)$.

Soit $\xi$ une section de $(T X_\R)\left| \right. _{\cup_{i=1}^j \lambda_i}$ de classe 
$C^1$, ne s'annulant pas, et transverse \`a $T (\cup_{i=1}^j \lambda_i)
\subset (T X_\R)\left| \right. _{\cup_{i=1}^j \lambda_i}$. 
Soit $x$ un point en lequel $\xi$ est tangent \`a $\lambda_i$ ($i \in \{ 1, 
\dots , j \}$), et $S_x$ une sph\`ere de $X_\R$ centr\'ee en $x$ choisie 
suffisament petite de sorte qu'elle soit incluse dans la r\'eunion 
$B_{k_1} \cup B_{k_2}$ des parties de $X_\R$
se situant de part et d'autre de $\lambda_i$, qu'elle ne contienne
aucun autre point en lequel $\xi$ est tangent \`a $\lambda_i$, et qu'elle ne soit
intersect\'ee par
$\lambda_i$ qu'en deux points $x_1$ et $x_2$. Choisissons de plus $S_x$ de 
sorte que 
$\xi (x_1)$ et $\xi (x_2)$ soient tangents \`a cette sph\`ere. Si ces deux vecteurs
ne se prolongent pas en un champ sur la sph\`ere, on pose
$o_i (\xi) (x) = 0$. Sinon, ces deux vecteurs induisent une orientation sur
$S_x$ et par suite une orientation sur le disque centr\'e en $x$ et bord\'e par 
$S_x$. Posons alors $o_i (\xi) (x) = \epsilon_1
t_{k_1} + \epsilon_2 t_{k_2} + q \epsilon_3 $ o\`u $\epsilon_p = +1$ 
(resp. $-1$) si 
l'orientation de $B_{k_p}$ co\"{\i}ncide (resp. ne co\"{\i}ncide pas) avec 
l'orientation de ce disque ($p \in \{ 1,2\}$) et $\epsilon_3 = 0$ (resp. $1$)
si $\lambda_i$ n'est pas incluse (resp. est incluse) dans $\Lambda_\R$.
Notons $o_i (\xi) = \sum o_i (\xi) (x)$, la somme \'etant prise sur les points 
$x$ en lesquels $\xi$ est tangent 
\`a $\lambda_i$ et orient\'e dans le sens de $\lambda_i$ ; et enfin
$i (\xi) = \sum_{i=1}^j s_i o_i (\xi)$.

Supposons qu'en les points d'intersections de $\cup_{i=1}^j \lambda_i$
et $\cup_{m=1}^n v^m$, le champ $\xi$ n'est pas tangent aux courbes $v^m$.
D'apr\`es le lemme 2.4, l'entier $i_\xi (v^m_\R , {\cal A}_{q, l})$ ne
d\'epend pas du prolongement de $\xi$ \`a $v^m_\R$. Notons $i_\xi (
{\cal A}_{q, l}) \in \Z^n$ le vecteur dont les coordonn\'ees sont
$\frac{1}{2} (i_\xi (v^m_\R , {\cal A}_{q, l}) + v^m \circ q \Lambda)$.

Remarquons enfin qu'un flot associ\'e au champ de vecteur $-J (\xi)$ permet 
de d\'ecaler les courbes $\lambda_i$ en des courbes $\tilde{\lambda}_i$ incluses
dans $X \setminus (X_\R \cup \cup_{m=1}^n v^m )$.

\begin{theo}
\label{multiple}
Supposons que $X$ satisfasse aux hypoth\`eses du th\'eor\`eme 2.5 et soit $\Lambda$ 
une $J$-courbe r\'eelle s\'eparante 
de $X$ telle que $q [\Lambda_\R] = q l \sum_{i=1}^j s_i [\lambda_i]
\in H_1 (X_\R ; \Z)$. Supposons qu'il existe un champ de vecteurs $\xi$ 
d\'efini sur $\cup_{i=1}^j 
\lambda_i$ comme ci-dessus, tel que les courbes $\tilde{\lambda}_i$ associ\'ees
satisfassent $\sum_{i=1}^j s_i [\tilde{\lambda}_i] = 0 \in H_1 (X \setminus 
(X_\R \cup \cup_{m=1}^n v^m) ; \Z)$, alors :
$$\int_{X_\R} \ind^2 d\chi + q l i(\xi) = - q_c \big( 
 [Q^{-1} (i_\xi ({\cal A}_{q , l} ))] \big) \quad 
\mod (2 q^2 l)$$
\end{theo}

{\bf D\'emonstration :}

Consid\'erons une extension du champ $\xi$ (toujours not\'ee $\xi$) au cycle
${\cal A}_{q , l}$ tout entier, et d\'ecalons ce cycle \`a l'aide d'un flot
associ\'e au champ $-J (\xi)$ de fa\c{c}on \`a obtenir un cycle simplicial
$\widetilde{\cal A}_{q ,l}$ g\'en\'erique. Le bord de ce cycle 
$\widetilde{\cal A}_{q ,l}$ est $q l \sum_{i=1}^j s_i 
\tilde{\lambda}_i $, et par hypoth\`ese il existe une $2$-cha\^{\i}ne 
$\widetilde{S}$ de $X \setminus (X_\R \cup \cup_{m=1}^n v^m) $ telle que
$\partial \widetilde{S} = \sum_{i=1}^j s_i \tilde{\lambda}_i $. A l'aide du flot pr\'ec\'edent, on peut prolonger cette cha\^{\i}ne $\widetilde{S}$ en une
$2$-cha\^{\i}ne $S$ bord\'ee par $\sum_{i=1}^j s_i \lambda_i $. Le cycle 
${\cal A} = {\cal A}_{q ,l} - q l S$ est alors un $2$-cycle entier de
$X$ qui rel\`eve ${\cal A}_{q , l}$.

La classe d'homologie de ${\cal A}$ se calcule en prolongeant le champ $\xi$
en un champ sur $\cup_{m=1}^n v^m$ et en d\'ecalant ces courbes \`a l'aide d'un 
flot associ\'e au champ de vecteurs $J (\xi)$ de fa\c{c}on \`a obtenir des courbes 
$\tilde{v}^m$
transverses \`a ${\cal A}$. Ces courbes $\tilde{v}^m$ n'intersectent pas $S$
par construction, et l'indice d'intersection $\tilde{v}^m \circ {\cal A}$ 
se calcule comme dans la d\'emonstration du th\'eor\`eme 2.5 et vaut
$\frac{1}{2} (v^m \circ q \Lambda + i_\xi (v^m_\R , 
{\cal A}_{q l}))$.
Il en d\'ecoule que $[{\cal A}] = \big( [Q^{-1} 
(i_\xi ({\cal A}_l ))] \big) \in H_2 (X ; \Z)$.

Consid\'erons \`a pr\'esent un champ continu $\zeta$ de vecteurs tangents \`a
$\Lambda_\R \cup \cup_{i=1}^j \lambda_i$, compatible avec l'orientation de
cette courbe, et prolongeons ce champ \`a ${\cal A}$ tout entier. D\'ecalons alors
le cycle ${\cal A}$ \`a l'aide d'un flot associ\'e au champ
$J (\zeta)$ de fa\c{c}on \`a obtenir un cycle $\widetilde{\cal A}$ imerg\'e et 
g\'en\'erique. L'intersection $\widetilde{\cal A} \circ c(\widetilde{\cal A})$
se calcule g\'eom\'etriquement et vaut $- \int_{X_\R} \ind^2 d\chi - q l 
i(\xi) \quad \mod (2 q^2 l)$. En effet le premier terme provient des z\'eros 
du champ de vecteurs
$\zeta$, le second des points de $\cup_{i=1}^j \lambda_i$ en 
lesquels $\xi$ et $\zeta$ sont
positivement colin\'eaires. Les intersections entre $q l S$ et $q 
\Lambda \cup q l c(S)$ qui sont \'eloign\'ees de $X_\R$ n'interviennent pas 
puisqu'elles sont coupl\'ees par $c$, donc nulles modulo
$2 q^2 l$. Le second terme peut s'obtenir en calculant
le nombre d'entrelacement de l'intersection des cycles $\widetilde{\cal A}$
et $c(\widetilde{\cal A})$ avec une petite sph\`ere centr\'ee en un tel point.
Cet entrelac est repr\'esent\'e par la figure suivante :
$$\vcenter
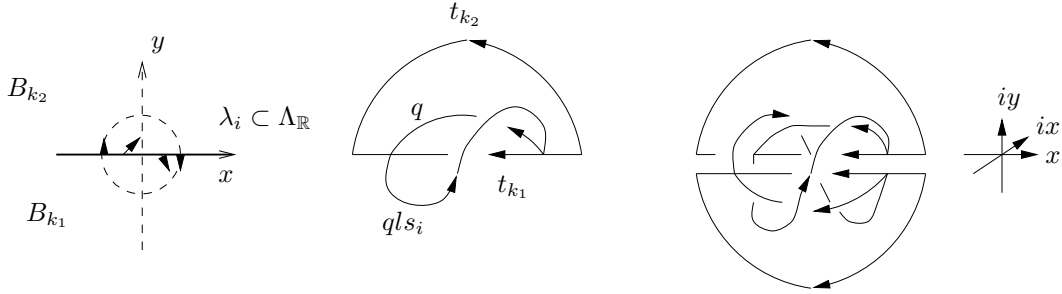
{\hbox{\begin{picture}(0,0)%
\epsfig{file=regl0.pstex}%
\end{picture}%
\setlength{\unitlength}{0.00066700in}%
\begingroup\makeatletter\ifx\SetFigFont\undefined%
\gdef\SetFigFont#1#2#3#4#5{%
  \reset@font\fontsize{#1}{#2pt}%
  \fontfamily{#3}\fontseries{#4}\fontshape{#5}%
  \selectfont}%
\fi\endgroup%
\begin{picture}(8112,2264)(1201,-5619)
\put(4651,-3511){\makebox(0,0)[lb]{\smash{\SetFigFont{10}{12.0}{\rmdefault}{\mddefault}{\updefault}$t_{k_2}$}}}
\put(5026,-4861){\makebox(0,0)[lb]{\smash{\SetFigFont{10}{12.0}{\rmdefault}{\mddefault}{\updefault}$t_{k_1}$}}}
\put(4126,-5161){\makebox(0,0)[lb]{\smash{\SetFigFont{10}{12.0}{\rmdefault}{\mddefault}{\updefault}$qls_i$}}}
\put(4351,-4261){\makebox(0,0)[lb]{\smash{\SetFigFont{10}{12.0}{\rmdefault}{\mddefault}{\updefault}$q$}}}
\put(2326,-3736){\makebox(0,0)[lb]{\smash{\SetFigFont{10}{12.0}{\rmdefault}{\mddefault}{\updefault}$y$}}}
\put(2851,-4786){\makebox(0,0)[lb]{\smash{\SetFigFont{10}{12.0}{\rmdefault}{\mddefault}{\updefault}$x$}}}
\put(1351,-5086){\makebox(0,0)[lb]{\smash{\SetFigFont{10}{12.0}{\rmdefault}{\mddefault}{\updefault}$B_{k_1}$}}}
\put(1201,-4111){\makebox(0,0)[lb]{\smash{\SetFigFont{10}{12.0}{\rmdefault}{\mddefault}{\updefault}$B_{k_2}$}}}
\put(2851,-4336){\makebox(0,0)[lb]{\smash{\SetFigFont{10}{12.0}{\rmdefault}{\mddefault}{\updefault}$\lambda_i  \subset \Lambda_\R$}}}
\put(8926,-4186){\makebox(0,0)[lb]{\smash{\SetFigFont{10}{12.0}{\rmdefault}{\mddefault}{\updefault}$iy$}}}
\put(9301,-4636){\makebox(0,0)[lb]{\smash{\SetFigFont{10}{12.0}{\rmdefault}{\mddefault}{\updefault}$x$}}}
\put(9226,-4411){\makebox(0,0)[lb]{\smash{\SetFigFont{10}{12.0}{\rmdefault}{\mddefault}{\updefault}$ix$}}}
\end{picture}
}}$$
(La sph\`ere $S^3$ priv\'ee de son point d'ordonn\'ee $y$ maximale a \'et\'e projet\'ee
 sur $\R^3$). $\square$

\subsubsection{Cas de deux courbes s\'eparantes disjointes.}

Soit $X$ une vari\'et\'e presque-complexe compacte de dimension $4$ (r\'eelle)
munie d'une structure r\'eelle $c$ et de l'orientation induite par la structure 
presque-complexe. Choisissons des $J$-courbes r\'eelles $v^1, \dots , 
v^n$ dans $X$. 

Soient $\Lambda$ et $e$ deux $J$-courbes r\'eelles s\'eparantes de $X$, 
transverses aux courbes $v^1, \dots , v^n$, et 
$\Lambda^+$ (resp. $e^+$) une partie de $\Lambda \setminus \Lambda_\R$ 
(resp. $e \setminus e_\R$) bord\'ee par $\Lambda_\R$ (resp. $e_\R$), qui induit
une orientation complexe sur $\Lambda_\R$ (resp. $e_\R$) que nous fixons.
Supposons que $e \circ \Lambda = 0$ et qu'il existe deux entiers $q,l$ 
tels que $q [\Lambda_\R] = l[e_\R] \in H_1 (X ; \Z)$. Notons $(B_k)$ les 
parties du compl\'ementaire de 
$\Lambda_\R \cup e_\R$ dans $X_\R$, que l'on oriente lorsqu'elles sont 
orientables, et choisissons une famille d'entiers $(t_k)$ tels que $t_k = 0$ 
si 
$B_k$ n'est pas orientable, et $\partial \sum_k t_k B_k = -q\Lambda_\R - l
e_\R$.
De la m\^eme fa\c{c}on que dans la section $2.1$, notons $\ind$ la fonction associ\'ee
\`a cette famille d'entiers, et ${\cal A}$ le $2$-cycle entier $q\Lambda^+ + 
le^+ + \sum_k t_k B_k$.

Notons enfin $i({\cal A})$ l'\'el\'ement de $\Z^n$ dont les coordonn\'ees sont
$\frac{1}{2} (q v^m \circ \Lambda + l v^m \circ e + i(v^m_\R ,{\cal A}))$
($m \in \{ 1, \dots , n \} $), et $Q$ la matrice de la forme d'intersection de
Poincar\'e dans la famille $([v^1], \dots , [v^n])$. (L'indice 
$i(v^m_\R ,{\cal A})$ a \'et\'e d\'efini dans la section $2.2.1$. Ici, la partie 
$J$-holomorphe de ${\cal A}$ est $q\Lambda^+ + le^+$, et $i_\xi (v^m_\R ,
{\cal A})(x)$ vaut $\pm q $ ou $\pm l$ selon que $x$ appartient \`a $v^m_\R \cap
\Lambda_\R$ ou $v^m_\R \cap e_\R$. Par ailleurs, $\cup_{i=1}^j \lambda_i$ est
vide, donc d'apr\`es le lemme $2.4$, l'entier $i_\xi (v^m_\R , {\cal A})$ est
ind\'ependant du champ $\xi$.)

\begin{theo}
Supposons que  $([v^1], \dots , [v^n])$ forme une base de l'espace 
$H_2 (X ; \Z)$, que $e \circ \Lambda = 0$ et qu'il existe deux entiers $q,l$ 
tels que $q [\Lambda_\R] = l[e_\R] \in H_1 (X ; \Z)$. En choisissant une
fonction $\ind$ comme ci-dessus et en notant ${\cal A}$ le $2$-cycle
associ\'e, on a :
$$\int_{X_\R} \ind^2 d\chi = -q_c \big( [Q^{-1} 
(i ({\cal A} ))]  \big).$$
\end{theo}

{\bf D\'emonstration :}

Consid\'erons un champ de vecteurs $\xi$ tangents \`a $\Lambda_\R \cup e_\R$,
compatible avec l'orientation de cette courbe, qu'on \'etend au cycle ${\cal A}$
tout entier. D\'ecalons le cycle ${\cal A}$ \`a l'aide d'un flot associ\'e au champ
$J (\xi)$ de fa\c{c}on \`a obtenir un cycle simplicial $\widetilde{\cal A}$ 
homologue \`a ${\cal A}$ et g\'en\'erique (v. \S $1.2$ pour une d\'efinition). De la 
m\^eme fa\c{c}on que dans $\cite{Rokh}$, on s'aper\c{c}oit que
${\cal A} \circ c({\cal A}) = - \int_{X_\R} \ind^2 d\chi$.

Prolongeons \`a pr\'esent le champ $\xi$ en un champ sur les courbes $v^m$, 
et d\'ecalons ces courbes \`a l'aide d'un flot associ\'e au champ
$J (\xi)$ de fa\c{c}on \`a obtenir des courbes $\tilde{v}^m$ homologues aux courbes
$v^m$ et transverses \`a ${\cal A}$. On s'aper\c{c}oit que $v^m \circ {\cal A} =
\frac{1}{2} (q v^m \circ \Lambda + l v^m \circ e + i(v^m_\R ,{\cal A}))$ 
et il en d\'ecoule que $[{\cal A}] = [Q^{-1} (i({\cal A}))] \in H_2 (X ; \Z)$. 
$\square$

\subsection{Quelques cons\'equences.}

\subsubsection{Formule de Rokhlin.}

Dans le cas o\`u l'entier $l$ est nul, les congruences donn\'ees par les 
th\'eor\`emes 2.2, 2.5, 2.8 et le corollaire 2.6 sont en fait des \'egalit\'es 
entre entiers. La formule donn\'ee par le
th\'eor\`eme 2.2 est alors la formule de Zvonilov (v. \cite{Zvo}). Si $X= \C P^2$,
la formule donn\'ee par le corollaire 2.6 est la formule des
orientations complexes de Rokhlin (v. \cite{Rokh}), et celle donn\'ee par 
le th\'eor\`eme \ref{multiple} est la formule de Mishachev (v. \cite{Mish}).
Remarquons \'egalement que lorsque $X$ est la surface r\'egl\'ee $\Sigma_1$ \`a
base $\cp$ (c'est le plan projectif \'eclat\'e en un point), la formule donn\'ee
par le th\'eor\`eme $2.9$ est la formule de Rokhlin.

\subsubsection{Congruences d'Arnol'd et de Mikhalkin.}

Si $X = \C P^2$, le corollaire 2.7 est la congruence d'Arnol'd 
(v. \cite{Arn}), et dans 
le cas o\`u $X = \cp \times \cp$ muni du produit des conjugaisons complexes, ce
corrolaire ainsi que le corollaire 2.6
\'etendent des congruences obtenues par Mikhalkin 
(v. \cite{Mikh} th\'eor\`eme 5.5, 5.7, 5.8 et 5.9 et 6.1d).

\section{Cas de surfaces fibr\'ees.}

Soient $\Sigma$ une vari\'et\'e presque-complexe de dimension $4$ r\'eelle et
$\Delta$ une surface de Riemann r\'eelle s\'eparante. Supposons qu'il existe
une submersion $J$-holomorphe $\pi : \Sigma \to \Delta$, ainsi qu'une 
structure r\'eelle $c_\Sigma$ sur $\Sigma$ telle
que $\pi \circ c_\Sigma = c_\Delta \circ \pi$. La vari\'et\'e $\Sigma$ est alors
dite {\it fibr\'ee} de base $\Delta$, et chaque fibre est une
$J$-courbe de $\Sigma$. Supposons que soit la base $\Delta$,
soit les fibres de $\pi$ sont simplement connexes, que ce fibr\'e poss\`ede une
$J$-section $e$ stable par $c_\Sigma$, et fixons $([v], [e])$ comme base de 
l'espace $H_2 (\Sigma ; \Z) \cong \Z^2$, o\`u $[v]$ est la classe d'homologie 
d'une fibre $v$. Remarquons que les composantes connexes de $\Sigma_\R$ sont 
des tores ou des bouteilles de Klein et qu'au-dessus de chaque composante  
de $\Delta_\R$, il y a au moins une composante de $\Sigma_\R$.

Soit $\Lambda$ une $J$-courbe r\'eelle s\'eparante (pas n\'ecessairement connexe) de 
$\Sigma$, et $a,b$ 
deux entiers tels que $[\Lambda] = (a,b) \in H_2 (\Sigma ; \Z)$. Munissons
$\Lambda_\R$ d'une orientation complexe, et supposons que dans chaque 
composante $\Sigma^\alpha_\R$ de $\Sigma_\R$, il existe une courbe connexe 
orient\'ee $\lambda_\alpha$ soit
disjointe de $\Lambda_\R$, soit incluse dans $\Lambda_\R$ (auquel cas 
l'orientation de $\lambda_\alpha$ et de $\Lambda_\R$ co\"{\i}ncident), et un
entier $l_\alpha$ tel que $[\Lambda_\R \cap \Sigma^\alpha_\R] = l_\alpha [
\lambda_\alpha] \in H_1 (\Sigma^\alpha_\R ; \Z)$. Sous cette hypoth\`ese,
la courbe $\Lambda$ est dite {\it admissible}.\\

Si $C$ est une courbe orient\'ee de $\Sigma_\R$ au-dessus 
d'une composante connexe (orient\'ee) $D$ de $\Delta_\R$, la 
classe d'homologie r\'ealis\'ee par le projet\'e de $C$ dans $H_1 (D ; \Z) \cong \Z$
est appel\'ee {\it nombre d'enroulement} de $C$ et not\'ee $n(C)$.

\begin{lemme}
Une courbe orient\'ee lisse d'une bouteille de Klein contient au plus deux
composantes de nombre d'enroulement $\pm 1$ et ne contient pas de
composantes de nombre d'enroulement de module strictement plus grand que $2$.
\end{lemme}

{\bf D\'emonstration :}

Le rev\^etement double d'orientation de la bouteille de Klein est un tore, 
et la courbe se 
rel\`eve en une courbe lisse dans ce tore, orient\'ee et invariante par 
l'involution du rev\^etement.
Chaque composante non contractile de cette courbe est une courbe connexe 
simple orient\'ee qui est soit disjointe, soit confondue avec son image par
l'involution du rev\^etement ; elle est donc  
n\'ec\'essairement 
homologue -au signe pr\`es- au relev\'e d'une section de la bouteille de Klein, ou
bien \`a une fibre du tore. Les composantes qui sont coupl\'ees 
par deux
avec l'involution du rev\^etement se quotientent en les composantes de nombre 
d'enroulement $\pm 2$ ou  nul, et celles qui sont stables par l'involution se 
quotientent en les
composantes de nombre d'enroulement $\pm 1$. Il ne peut pas y avoir plus que 
deux telles composantes, puisque les voisinages tubulaires de ces composantes 
dans la bouteille de Klein sont des rubans de M\oe bius. $\square$\\

La condition d'existence de courbes $\lambda_\alpha$ \'equivaut \`a imposer
que dans chacune des composantes $\Sigma^\alpha_\R$ qui est une bouteille de 
Klein, la courbe $\Lambda_\R$ a au plus une composante de nombre d'enroulement
 $\pm 1$.

Notons $l = pgcd (l_\alpha)$,  
choisissons une famille $(t_k)$ d'entiers telle que si une courbe 
$\lambda_\alpha$ est bilat\`ere dans $\Sigma^\alpha_\R$, alors $t_k$ est nul 
d'un côt\'e de $\lambda_\alpha$, et notons $\ind$ la fonction associ\'ee \`a ce 
choix d'entiers (v. \S$2.1$ pour une d\'efinition).

Notons $s_\alpha$ l'entier tel que $[\Lambda_\R \cap \Sigma^\alpha_\R] = 
l s_\alpha [\lambda_\alpha] \in H_1 (\Sigma^\alpha_\R ; \Z)$ et 
$s = e_\R \circ \sum s_\alpha [\lambda_\alpha] \in \Z /2 \Z$ (la somme \'etant
prise sur tous les indices $\alpha$, ou de fa\c{c}on \'equivalente sur les indices 
$\alpha$ tels que $e_\R$ rencontre $\Sigma_\R^\alpha$.
 De m\^eme, notons $l t \in \Z$ la somme des nombres
d'enroulements des courbes $\Lambda_\R \cap \Sigma^\alpha_\R$, o\`u
$\Sigma^\alpha_\R$ est une composante rencontr\'ee par $v_\R$ (ce nombre
d'enroulement est bien divisible par $l$, ce qui justifie la notation).

\begin{theo}
Soit $\Sigma$ une vari\'et\'e fibr\'ee de dimension $4$ r\'eelle 
qui admette une $J$-section $e$ r\'eelle. Supposons que soit la base 
$\Delta$, soit les fibres de $\Sigma$ sont simplement connexes. 
Soit $\Lambda$ une $J$-courbe r\'eelle s\'eparante admissible de $\Sigma$, de 
classe $a[v] + b[e]$ dans 
$H_2 (\Sigma ; \Z)$ o\`u $v$ est la classe d'homologie d'une fibre. 
Choisissons une fonction $\ind$ (associ\'ee \`a un choix d'entiers $(t_k)$),
des entiers $l$, $t$ et un \'el\'ement $s \in \Z / 2\Z$ comme ci-dessus, et
notons $\nu = e \circ e$. On a :
$$\int_{\Sigma_\R} \ind^2 d \chi = \frac{1}{2} (b + lt)(a + ls + 
\frac{\nu}{2} (b - lt)) \quad \left[ 
\begin{array}{l} 
\mod (2l) \mbox{ si $l$ est pair,} \\
\mod (l) \mbox{ sinon.}
\end{array} \right.$$
En particulier, si $l$ est pair et $B$ est la r\'eunion des parties $B_k$
pour lesquelles $t_k$ est impair, 
$$\chi (B) = \frac{1}{2} (b + lt)(a + ls + 
\frac{\nu}{2} (b - lt)) \quad mod(4)$$
\end{theo}

{\bf D\'emonstration :}

Le corollaire $2.6$ s'applique, $Q =
\left(
\begin{array}{cc}
0 &1 \\
1 & \nu
\end{array}
\right)$ et
$i({\cal A}_{l}) = \frac{1}{2} (b + lt , a + \nu b + ls) \; mod(l)$. 
En effet, $i({\cal A}_{l})$ ne d\'epend que de ${\cal A}_{l} \cap \Sigma_\R$
et est invariant par isotopie de $\Sigma_\R$ (ce qui d\'ecoule 
du lemme $2.4$). En choisissant une isotopie de fa\c{c}on \`a rendre
${\cal A}_{l} \cap \Sigma_\R$ voisin de $\cup \lambda_\alpha$, on
obtient le r\'esultat. $\square$

\begin{theo}
\label{2comp}
Soit $\Sigma$ une surface r\'egl\'ee, de base s\'eparante, telle que $\Sigma_\R$ 
est une bouteille de Klein (en particulier, $\Sigma_\R$ est connexe) et
admettant une $J$-section $e$ r\'eelle. 
Soit $\Lambda$ une $J$-courbe r\'eelle s\'eparante de $\Sigma$ telle que
$\Lambda_\R$ contienne deux composantes de nombre d'enroulement $\pm 1$.
Notons $\lambda$ une de ces deux composantes et $l= n(\Lambda_\R) n(\lambda)$. 
Soient $a,b$ deux entiers et $s$ un \'el\'ement de $\Z / 2 \Z$ tels que 
$[\Lambda] = (a,b)
\in H_2 (\Sigma ; \Z)$ et $[\lambda] = (1+s)[v_\R] + n(\lambda) [e_\R] \in H_1 
(\Sigma_\R ; \Z)$ o\`u $v$ est une fibre r\'eelle. 
Soit $(t_k)_{k \in K }$ une famille d'entiers 
telle que $\partial (\sum_{k \in K} t_k B_k) = 2l\lambda - 2 \Lambda_\R$, 
en notant $\ind$ la fonction de $\Sigma_\R$
associ\'ee \`a cette famille (v. \S\ref{tk} pour une d\'efinition) et 
$\nu = e \circ e$, on a :
$$\int_{\Sigma_\R} ind^2 d\chi + 4l \max(0 ,n(\lambda)) = 2 (b+l)
(a +ls +\frac{\nu}{2}(b-l) ) \, \mod (8l)$$

\end{theo}

{\bf D\'emonstration :}

Il existe un champ $\xi$ de vecteurs tangents de $\Sigma_\R$ de classe 
$C^1$, d\'efini sur $\lambda$, transverse aux fibres, tel que projet\'e sur la base
$\Delta_\R$ du fibr\'e il induise un champ compatible avec l'orientation de 
$\Delta_\R$.
Soit $\tilde{\lambda}$ une courbe obtenue en d\'ecalant $\lambda$ \`a l'aide
d'un flot associ\'e au champ $-J \xi$. Cette courbe est plong\'ee dans la
restriction du fibr\'e au dessus d'une partie $\Delta^-$ de $\Delta \setminus 
\Delta_\R$ et n'intersecte pas $e$. Par suite, 
$[\tilde{\lambda}] = 0 \in H_1 (\Sigma \setminus (\Sigma_\R
\cup (e \cup v)) ; \Z)$, et le th\'eor\`eme \ref{multiple}
s'applique avec $q = 2$.

Or $i(\xi) = 2\max (0,n(\lambda))  \mod (4)$, 
$i_\xi ({\cal A}_{2, l}) = (b+l,a + \nu b +ls)  
\mod (4l)$, et $\mbox{$Q =
\left(
\begin{array}{cc}
0 &1 \\
1 & \nu
\end{array}
\right)$. $\square$}$

\begin{theo}
Soit $\Sigma_\nu$, $\nu > 0$, une surface r\'egl\'ee de base $\cp$, et $e$ son
diviseur exceptionnel ($e \circ e = - \nu$). Soit $\Lambda$ 
une $J$-courbe r\'eelle s\'eparante de classe $a [v] + b[e]$ dans 
$H_2 (\Sigma_\nu ; \Z)$ ($[v]$ est la classe d'homologie d'une fibre), 
qui n'intersecte pas le diviseur
exceptionel $e$ de $\Sigma_\nu$, et telle que $a$ et $b$ sont pairs. 
Il existe un entier
$l$ tel que $[\Lambda_\R] = l [e_\R] \in H_1 ((\Sigma_\nu)_\R ; \Z)$ (
$\Lambda_\R$ et $e_\R$ sont munies d'orientations complexes) et en choisissant
une fonction $\ind$ comme dans le $\S 2.3.2$, on a :
$$\int_{(\Sigma_\nu)_\R} \ind^2 d \chi = \frac{\nu}{4} (b^2 - l^2)$$
\end{theo}

{\bf D\'emonstration :}

Le th\'eor\`eme $2.9$ s'applique avec $a = \nu b$ et $q=1$. Choisissons 
\`a pr\'esent $([v] , 
[\tilde{e} ])$ comme base de l'espace $H_2 (\Sigma_\nu ; \Z)$, o\`u $v$ est une
 fibre au-dessus d'un point r\'eel, et $\tilde{e}$ une section holomorphe 
invariante par la conjugaison complexe et de carr\'e d'intersection $\nu$.
Dans cette base $Q =
\left(
\begin{array}{cc}
0 &1 \\
1 & \nu
\end{array}
\right)$, et $i({\cal A}) = \frac{1}{2} (b + l , a) $. $\square$

\section{Exemples.}

Dans cette section sont donn\'es quelques exemples d'application des r\'esultats
pr\'ec\'edents dans le cas de surfaces fibr\'ees de partie r\'eelle connexe.

\subsection{Lorsque $\Sigma_\R$ est orientable.}

\begin{prop}
Soit $\Sigma$ une vari\'et\'e satisfaisant aux hypoth\`eses du th\'eor\`eme $3.2$,
 telle que $\Sigma_\R$ est un tore.
Soit $\Lambda$ une $(M-2r)$ $J$-courbe r\'eelle 
de classe $(a,4)$ dans $\Sigma$, telle que $\Lambda_\R$ 
contienne deux composantes
 de classe $\pm(s,1)$. Soient $p$, $q$ les
nombres d'ovales se situant dans les deux composantes de $\Sigma_\R$ 
priv\'e de ces composantes de classe $\pm(s,1)$ et $\nu = e \circ e$. 
\begin{enumerate}
\item Si $2s = a - \nu \, \mod (4)$, et $p$, $q$ sont de parit\'e 
oppos\'ee \`a $\frac{a}{2} + r$, la courbe n'est pas s\'eparante.
\item Si la courbe est s\'eparante et $2s = a - \nu + 2 \, \mod (4)$,
alors $[\Lambda_\R] = 0 \in H_1 (\Sigma_\R ; \Z)$ si $p$, $q$ sont de 
m\^eme parit\'e que $\frac{a}{2} + r$, et $[\Lambda_\R] = \pm 2 (s,1) \in 
H_1 (\Sigma_\R ; \Z)$ sinon.
\end{enumerate}
\end{prop} 

{\bf D\'emonstration :}

Si $\Lambda$ est s\'eparante, le th\'eor\`eme $3.2$ s'applique, $b=4$, $t=1$, et 
$\chi (B)$ vaut le nombre
d'ovales plus $2q$ modulo $4$. Comme $\Lambda$ est une $(M-2r)$-courbe dont
la partie r\'eelle poss\`ede deux composantes non contractiles, le nombre d'ovales
vaut $g(\Lambda) + 1 - 2r - 2 = g(\Lambda) - 1 - 2r$. Le genre de $\Lambda$
se calcule \`a l'aide de la formule d'adjonction et vaut $(a-1)(b-1) +
\frac{1}{2} \nu b(b-1) + g(\Delta ) b$ avec $b=4$ ($\Delta$ est la base du 
fibr\'e). En appliquant le th\'eor\`eme $3.2$, on en d\'eduit que
$3(a-1) +6 \nu -1 -2r +2q = \frac{1}{2}(al + l^2 (s -\frac{\nu}{2})) \mod (4)$
(en remarquant que $\nu$ est pair). Il en d\'ecoule que
$q = \frac{a}{2} + r \mod (2)$ si $l=0$, et $q = s + r +\frac{\nu}{2}
\mod (2)$ si $l= \pm 2$. $\square$

\begin{prop}
Soit $\Sigma$ une vari\'et\'e satisfaisant aux hypoth\`eses du th\'eor\`eme $3.2$,
 telle que $\Sigma_\R$ est un tore.
Soit $\Lambda$ une $(M-2r)$ $J$-courbe r\'eelle s\'eparante
de classe $(a,6)$ dans $\Sigma$, telle que $\Lambda_\R$ 
contienne deux composantes
 de classe $\pm(s,2)$. Soient $p$, $q$ les
nombres d'ovales se situant dans les deux composantes de $\Sigma_\R$ 
priv\'e de ces composantes de classe $\pm(s,2)$ et $\nu = e \circ e$. Alors,
avec les notations du th\'eor\`eme $3.2$, 
 $[\Lambda_\R] = 0 \in 
H_1 (\Sigma_\R ; \Z)$ si $p$, $q$ sont de 
m\^eme parit\'e que $r+ \frac{1}{2} \chi (\Delta)$, et  $[\Lambda_\R] = \pm 2 
(s,2) \in H_1 (\Sigma_\R ; \Z)$ sinon.
\end{prop}

{\bf D\'emonstration :}

Le th\'eor\`eme $3.2$ s'applique, $b=6$, $t=2$, et $\chi (B)$ vaut le nombre
d'ovales plus $2q$ modulo $4$. Comme $\Lambda$ est une $(M-2r)$-courbe dont
la partie r\'eelle poss\`ede deux composantes non contractiles, le nombre d'ovales
vaut $g(\Lambda) + 1 - 2r - 2 = g(\Lambda) - 1 - 2r$. Le genre de $\Lambda$
se calcule \`a l'aide de la formule d'adjonction et vaut $(a-1)(b-1) +
\frac{1}{2} \nu b(b-1) + g(\Delta ) b$ avec $b=6$ ($\Delta$ est la base du 
fibr\'e). En appliquant le th\'eor\`eme $3.2$, on en d\'eduit que
$5(a-1) + 15 \nu +2g(\Delta ) -1 -2r +2q = 3a + 3ls +9\nu  \mod (4)$
(en remarquant que $\nu$ est pair), ou encore
$q = 1 - g(\Delta ) + r - \frac{1}{2} ls \mod (2)$. $\square$

\begin{prop}
Il n'existe pas de $(M-4)$-courbe alg\'ebrique r\'eelle s\'eparante de bidegr\'e 
$(8,8)$ sur l'hyperbolo\"{\i}de dont la partie r\'eelle contient $6$
composantes de classe $(0,1)$ et qui r\'ealise la classe d'isotopie de courbes 
suivante :
$$\vcenter
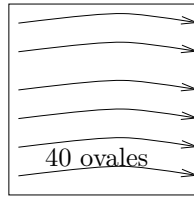
{\hbox{\begin{picture}(0,0)%
\epsfig{file=regle7.pstex}%
\end{picture}%
\setlength{\unitlength}{0.00066700in}%
\begingroup\makeatletter\ifx\SetFigFont\undefined%
\gdef\SetFigFont#1#2#3#4#5{%
  \reset@font\fontsize{#1}{#2pt}%
  \fontfamily{#3}\fontseries{#4}\fontshape{#5}%
  \selectfont}%
\fi\endgroup%
\begin{picture}(1524,1524)(4789,-3973)
\put(5101,-3736){\makebox(0,0)[lb]{\smash{\SetFigFont{10}{12.0}{\rmdefault}{\mddefault}{\updefault}$40$ ovales}}}
\end{picture}
}}$$
\end{prop}

{\bf D\'emonstration :}

Le th\'eor\`eme $3.2$ s'applique, $a=b=8$, $t=1$, $s=0$, $l=6$ et $\nu = 0$.
Or dans la congruence modulo $12$ donn\'ee par le th\'eor\`eme $3.2$, le membre
de gauche vaut $4$, et le membre de droite vaut $8$. $\square$\\

{\bf Remarques :}
Dans la proposition 4.1, si $\Sigma$ est un hyperbolo\"{\i}de, la partie 
$1$ se d\'eduit des 
r\'esultats ant\'erieurs de Mikhalkin (v. \cite{Mikh}). De m\^eme
dans la partie $2$, le fait que $[\Lambda_\R] = 0 \in H_1 (\Sigma_\R ; 
\Z)$ implique que $p$, $q$ sont de m\^eme parit\'e que $\frac{a}{2} - r$ peut se 
d\'eduire de la formule de Zvonilov (v. \cite{Zvo}). 

Remarquons qu'il est facile d'\'enoncer une proposition analogue pour les
courbes de classe $(a,b)$ avec $a$, $b$ pairs, dont la partie r\'eelle contient
$(b-2)$ composantes de classe $\pm(s,1)$.

Dans la proposition 4.2, le fait que $[\Lambda_\R] = 0 \in H_1 
(\Sigma_\R ; \Z)$ implique que $p$, $q$ sont de parit\'e oppos\'ee \`a
$r+ \frac{1}{2} \chi (\Delta)$ peut se d\'eduire de la formule de Zvonilov 
(v. \cite{Zvo}).

\subsection{Lorsque $\Sigma_\R$ n'est pas orientable.}

\begin{prop}
Soit $\Sigma$ une vari\'et\'e satisfaisant aux hypoth\`eses du th\'eor\`eme $3.2$,
 telle que $\Sigma_\R$ est une bouteille de Klein.
Soit $\Lambda$ une $(M-2r)$ $J$-courbe r\'eelle 
de classe $(a,4)$ dans $\Sigma$, telle que $\Lambda_\R$ 
contienne  une composante de nombre d'enroulement $\pm 2$. La composante de
nombre d'enroulement $\pm 2$ s\'epare $\Sigma_\R$ en deux rubans de
M\oe bius, et l'\^ame de l'un d'entre eux a un nombre de points d'intersections 
impair avec $e_\R$. Notons $q$ le nombre d'ovale se situant dans ce ruban et
 $\nu = e \circ e$.
Si  $q = \frac{1}{2} (\nu - 1)  + r \mod (2)$, alors la courbe n'est pas 
s\'eparante. 
\end{prop} 

{\bf D\'emonstration :}

Si $\Lambda$ est s\'eparante, le th\'eor\`eme $3.2$ s'applique, $b=4$, $l=\pm 2$,
choisissons $\lambda$ de sorte que $s=1$ et 
$\chi (B)$ vaut le nombre
d'ovales plus $2q$ modulo $4$. Comme $\Lambda$ est une $(M-2r)$-courbe dont
la partie r\'eelle poss\`ede une composante non contractile, le nombre d'ovales
vaut $g(\Lambda) + 1 - 2r - 1 = g(\Lambda) - 2r$. Le genre de $\Lambda$
se calcule \`a l'aide de la formule d'adjonction et vaut $(a-1)(b-1) +
\frac{1}{2} \nu b(b-1) + g(\Delta ) b$ avec $b=4$ ($\Delta$ est la base du 
fibr\'e). En appliquant le th\'eor\`eme $3.2$, on en d\'eduit que
$3(a-1) +6 \nu -2r +2q = -a +2 - \nu \mod (4)$. Il en d\'ecoule que
$2q = 1 + \nu + 2r \mod (4)$. $\square$

\begin{prop}
Soit $\Sigma$ une vari\'et\'e satisfaisant aux hypoth\`eses du th\'eor\`eme $3.2$,
 telle que $\Sigma_\R$ est une bouteille de Klein.
Soit $\Lambda$ une $(M-2r)$ $J$-courbe r\'eelle 
de classe $(a,6)$ dans $\Sigma$, telle que $\Lambda_\R$ 
contienne deux composantes de nombre d'enroulement $\pm 2$. Ces composantes 
s\'eparent $\Sigma_\R$ en deux rubans de M\oe bius et un cylindre, notons
$q$ le nombre d'ovales se situant dans les rubans de  M\oe bius.
Si  $q = r - g(\Delta )  \mod (2)$, la courbe n'est pas s\'eparante. 
\end{prop}

{\bf D\'emonstration :}

Le th\'eor\`eme $3.2$ s'applique, $b=6$, $l=0$ ou $\pm 4$,
choisissons $\lambda$ de sorte que $s=1$ et 
$\chi (B)$ vaut le nombre
d'ovales plus $2p$ modulo $4$ (o\`u $p$ est le nombre d'ovales se situant 
dans le cylindre de $\Sigma_\R$ priv\'e des composantes non contractiles de
$\Lambda_\R$). Comme $\Lambda$ est une $(M-2r)$-courbe dont
la partie r\'eelle poss\`ede deux composantes non contractiles, le nombre d'ovales
vaut $g(\Lambda) + 1 - 2r - 2 = g(\Lambda) - 1 - 2r$. Le genre de $\Lambda$
se calcule \`a l'aide de la formule d'adjonction et vaut $(a-1)(b-1) +
\frac{1}{2} \nu b(b-1) + g(\Delta ) b$ avec $b=6$ ($\Delta$ est la base du 
fibr\'e). En appliquant le th\'eor\`eme $3.2$, on en d\'eduit que
$5(a-1) + 15 \nu +2g(\Delta ) -1 -2r +2p = -a + \nu  \mod (4)$, ou encore
$p = 1 + \nu - g(\Delta ) + r = r - g(\Delta ) \mod (2)$. $\square$

\begin{prop}
Il n'existe pas de courbe alg\'ebrique r\'eelle maximale de classe $(8,6)$
 sur la surface $\Sigma_1$ dont la partie r\'eelle contient deux 
composantes de nombre d'enroulement $2$ et qui r\'ealise la classe d'isotopie 
de courbes suivante :
$$\vcenter
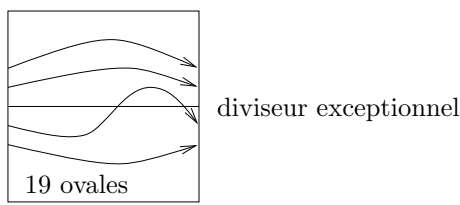
{\hbox{\begin{picture}(0,0)%
\epsfig{file=regle8.pstex}%
\end{picture}%
\setlength{\unitlength}{0.00066700in}%
\begingroup\makeatletter\ifx\SetFigFont\undefined%
\gdef\SetFigFont#1#2#3#4#5{%
  \reset@font\fontsize{#1}{#2pt}%
  \fontfamily{#3}\fontseries{#4}\fontshape{#5}%
  \selectfont}%
\fi\endgroup%
\begin{picture}(1662,1524)(3289,-4573)
\put(4951,-3886){\makebox(0,0)[lb]{\smash{\SetFigFont{10}{12.0}{\rmdefault}{\mddefault}{\updefault}diviseur exceptionnel}}}
\put(3451,-4486){\makebox(0,0)[lb]{\smash{\SetFigFont{10}{12.0}{\rmdefault}{\mddefault}{\updefault}$19$ ovales}}}
\end{picture}
}}$$
\end{prop}

{\bf D\'emonstration :}

Le th\'eor\`eme $3.2$ s'applique, $a=8$, $b=6$, $l=4$, choisissons $\lambda$ de 
sorte que $s=1$ et $\nu = -1$.
Or dans la congruence modulo $8$ donn\'ee par le th\'eor\`eme $3.2$, le membre
de gauche vaut $3$, et le membre de droite vaut $-1$. $\square$

\begin{prop}
Les sch\'emas complexes suivants ne sont pas r\'ealis\'es  sur la surface 
$\Sigma_1$ comme parties r\'eelles de
courbes alg\'ebrique r\'eelles maximales de classe $(5,4)$ contenant deux
composantes de nombre d'enroulement $1$ :
$$\vcenter{\hbox{\begin{picture}(0,0)%
\epsfig{file=regle9.pstex}%
\end{picture}%
\setlength{\unitlength}{0.00066700in}%
\begingroup\makeatletter\ifx\SetFigFont\undefined%
\gdef\SetFigFont#1#2#3#4#5{%
  \reset@font\fontsize{#1}{#2pt}%
  \fontfamily{#3}\fontseries{#4}\fontshape{#5}%
  \selectfont}%
\fi\endgroup%
\begin{picture}(5862,1524)(2089,-4873)
\put(7951,-4261){\makebox(0,0)[lb]{\smash{\SetFigFont{10}{12.0}{\rmdefault}{\mddefault}{\updefault}diviseur exceptionnel}}}
\end{picture}
}}$$
\end{prop}

{\bf D\'emonstration :}

Le th\'eor\`eme \ref{2comp} s'applique, $a=5$, $b=4$, $l=2$, choisissons $\lambda$
 de sorte que $s=1$ et $\nu = -1$. Or dans la congruence modulo $8$ donn\'ee 
par le th\'eor\`eme \ref{2comp}, pour ces courbes, le membre
de gauche vaut $8$, et le membre de droite vaut $0$. $\square$

\end{document}